\documentclass{article}

\newcommand{\be}{\begin{equation}}
\newcommand{\ee}{\end{equation}}
\newcommand{\bea}{\begin{eqnarray}}
\newcommand{\eea}{\end{eqnarray}}
\newcommand{\barray}{\begin{array}}
\newcommand{\earray}{\end{array}}
\newcommand{\pa}{\partial}
\newcommand{\nn}{\nonumber}
\newcommand{\bitem}{\begin{itemize}}
\newcommand{\eitem}{\end{itemize}}
\newtheorem{teo}{Theorem}[section]
\newcommand{\bt}{\begin{teo}}
\newcommand{\et}{\end{teo}}
\newtheorem{Def}{Definition}[section]
\newcommand{\bd}{\begin{Def}}
\newcommand{\ed}{\end{Def}}
\newtheorem{lem}{Lemma}[section]
\newcommand{\bl}{\begin{lem}}
\newcommand{\el}{\end{lem}}
\newtheorem{prop}{Proposition}[section]
\newcommand{\bp}{\begin{prop}}
\newcommand{\ep}{\end{prop}}
\newtheorem{cor}{Corollary}[section]
\newcommand{\bc}{\begin{cor}}
\newcommand{\ec}{\end{cor}}
\newtheorem{ex}{Example}[section]
\newcommand{\bex}{\begin{ex}}
\newcommand{\eex}{\end{ex}}
\newtheorem{rem}{Remark}[section]
\newcommand{\br}{\begin{rem}}
\newcommand{\er}{\end{rem}}

\catcode `\@=11
\@addtoreset{equation}{section}

\begin{document}

\begin{center}
{\Large \textbf{Quasi-Frobenius algebras and their
integrable \\ $N$-parametric deformations
generated by compatible \\ $(N \times N)$-metrics of
constant Riemannian curvature\footnote{This work was
supported by the Alexander von Humboldt Foundation
(Germany), the Russian Foundation for Basic Research
(Grant No. 02--01--00803), and INTAS (Grant No. 99--1782).}}}
\end{center}

\medskip

\centerline{\large {O. I. Mokhov}}
\bigskip

\section{Introduction}

In this paper, we prove that the description
of pencils of compatible $(N \times N)$-metrics
of constant Riemannian curvature is equivalent
to a special class of integrable
$N$-parametric deformations of quasi-Frobenius
(in general, noncommutative) algebras.

A finite-dimensional algebra ${\cal {Q}}$
is said to be
{\it quasi-Frobenius} if the identity
\be
(a b) c =
(a c) b,  \ \ \ a, b, c \in {\cal {Q}} \label{q1}
\ee
is fulfilled and also an invariant nondegenerate
symmetric bilinear form
$<a, b>$ is given,
\be
<a b, c> = <a, c b>, \ \ \ <a, b>=<b, a>,
\ \ \ a, b, c \in {\cal {Q}}. \label{q2}
\ee
Recall that a finite-dimensional commutative
associative algebra equipped with an
invariant nondegenerate symmetric bilinear form
is called {\it a Frobenius algebra}
(here, we do not require an existence
of a unit in Frobenius algebra).
Any commutative quasi-Frobenius algebra
is always Frobenius, i.e., if the identity
$a b =  b a$
({\it commutativity}) is fulfilled in a quasi-Frobenius
algebra, then the identities
\be
(a b) c =
a (b c) \ \ \  ({\it associativity}), \label{f1}
\ee
\be
<a b, c> = <a, b c> \ \ \  ({\it invariance
\ of \ bilinear \ form}) \label{f2}
\ee
are also always fulfilled in the algebra.

Identity (\ref{q1}) meaning commutativity
of the operators of right-sided multiplication
in the algebra: $R_a R_b = R_b R_a$, where $R_a b = b a$,
naturally arises as one of identities in algebras
describing some special classes of first-order Poisson
brackets linear in field variables
(see \cite{1x}--\cite{3x}).
The theory of left-symmetric algebras with
the additional identity (\ref{q1}), which correspond
to linear one-dimensional Poisson brackets of
hydrodynamic type and are called
{\it Novikov algebras}, was developed in
\cite{4x}--\cite{6x} ({\it left-symmetric algebras} or
{\it Vinberg algebras}, i.e.,
algebras with the identity $ a (b c) - (a b) c = b (a c) - (b a) c$
or, equivalently,
$L_a L_b - L_b L_a = L_{a b - ba}$, where $L_a b = a b,$
were considered in \cite{7x}).
Identity (\ref{q1}) naturally arises also in
algebras describing all multidimensional
Poisson brackets of hydrodynamic type (the algebras
were studied by the present author in
\cite{8x}).

Introducing the considered in this paper notion
of quasi-Frobenius algebras is motivated
by Dubrovin's theory of Frobenius manifolds
\cite{9x} and by a natural generalization of
Frobenius manifolds (quasi-Frobenius manifolds)
that is related to arbitrary flat pencils of metrics
\cite{10x}, see also \cite{11x}, \cite{12x}.
Dubrovin proved that two-dimensional topological
field theories (Frobenius manifolds)
correspond to a special class of integrable
deformations of Frobenius algebras, and also to
a special class of quasihomogeneous flat pencils of
metrics (see \cite{9x}, \cite{12xx}).
General flat pencils of metrics locally correspond to
quasi-Frobenius manifolds
(manifolds with quasi-Frobenius structure on
tangent spaces, see \cite{10x}) and, respectively,
to special deformations of quasi-Frobenius algebras.
Integrability of nonlinear equations
describing these deformations of quasi-Frobenius algebras
was proved by the present author in
\cite{13x}, \cite{14x} by the method of the inverse scattering
problem, see also \cite{11x}, \cite{15x}, \cite{16x}
(in \cite{21x}, a Lax pair for these equations
was also indicated). In the present work,
we prove that pencils of compatible metrics
of constant Riemannian curvature
also generate special integrable deformations of
quasi-Frobenius algebras.
We very hope that integrable deformations of
quasi-Frobenius algebras, constructed in the paper,
will also prove to be useful in two-dimensional
topological field theories.
In particular, we conjecture that the deformations
of noncommutative associative algebras, recently found
by S.M.Natanzon in \cite{22x} (see also \cite{23x})
and corresponding to
open-closed two-dimensional
topological field theories, are integrable
by the method of the inverse scattering problem
and are special reductions of the considered
class of integrable deformations
of noncommutative quasi-Frobenius algebras.
In this connection, the study of deformations of
noncommutative quasi-Frobenius algebras and, especially,
the extraction of associative deformations among them
is of a great interest.

\section{Compatible metrics of constant \\
Riemannian curvature}

Recall that two pseudo-Riemannian contravariant
metrics
$g_1^{ij} (u)$ and $g_2^{ij} (u)$ are said to be
{\it compatible} if for any linear combination
$g^{ij} (u) = \lambda_1 g_1^{ij} (u) + \lambda_2 g_2^{ij} (u)$
of these metrics, where $\lambda_1$ and $\lambda_2$
are arbitrary constants such that
$\det ( g^{ij} (u) ) \not\equiv 0$,
the coefficients of the corresponding Levi-Civita
connections and the components of the corresponding
Riemannian curvature tensors are related by the same
linear formula:
$\Gamma^{ij}_k (u) = \lambda_1 \Gamma^{ij}_{1, k} (u) +
\lambda_2 \Gamma^{ij}_{2, k} (u)$ and
$R^{ij}_{kl} (u) = \lambda_1 R^{ij}_{1, kl} (u)
+ \lambda_2 R^{ij}_{2, kl} (u)$
(in this case, we also say that {\it the metrics
$g^{ij}_1 (u)$ and $g^{ij}_2 (u)$ form a pencil
of compatible metrics}) \cite{16x}, \cite{24x}.
Flat pencils of metrics, which are none other than
compatible nondegenerate local Poisson brackets
of hydrodynamic type (compatible Dubrovin--Novikov
brackets \cite{25x}), were introduced in \cite{9x}.
Two pseudo-Riemannian contravariant metrics
$g_1^{ij} (u)$ and $g_2^{ij} (u)$ of constant
Riemannian curvatures
$K_1$ and $K_2$, respectively, are said to be {\it compatible}
if any linear combination
$g^{ij} (u) = \lambda_1 g_1^{ij} (u) + \lambda_2 g_2^{ij} (u)$
of these metrics,
where $\lambda_1$ and $\lambda_2$ are arbitrary constants
such that
$\det ( g^{ij} (u) ) \not\equiv 0$, is a metric
of constant Riemannian curvature
$\lambda_1 K_1 + \lambda_2 K_2$,
and the coefficients of the corresponding
Levi-Civita connections are related by the same
linear formula:
$\Gamma^{ij}_k (u) = \lambda_1 \Gamma^{ij}_{1, k} (u) +
\lambda_2 \Gamma^{ij}_{2, k} (u)$ \cite{16x}, \cite{24x}.
In this case, we also say that {\it the metrics
$g_1^{ij} (u)$ and $g_2^{ij} (u)$ form a pencil
of compatible metrics of constant Riemannian metrics}
\cite{16x}, \cite{24x}.
It is obvious that all these definitions are consistent
with one another: if compatible metrics are metrics of
constant Riemannian curvature, then they form a pencil
of compatible metrics of constant Riemannian curvature,
and if compatible metrics are flat, then they form
a flat pencil of metrics.

In \cite{26x} nonlocal Poisson brackets of
hydrodynamic type that have the following form
(the Mokhov--Ferapontov brackets):

\be
\{ I,J \} = \int {\delta I \over \delta u^i(x)}
\left ( g^{ij}(u(x)) {d \over dx} + b^{ij}_k (u(x))\, u^k_x
+ K u^i_x \left ( {d \over dx} \right )^{-1} u^j_x \right )
{\delta J \over \delta u^j(x)} dx,
\label{nonl}
\ee
where $I [u]$ and $J [u]$ are arbitrary functionals
on the space of functions (fields) $u^i (x),$ $1 \leq i \leq N,$
of one independent variable $x$, $u = (u^1,..., u^N)$ are
local coordinates on a certain given smooth
$N$-dimensional manifold $M$, the coefficients $g^{ij} (u)$ and
$b^{ij}_k (u)$ of bracket (\ref{nonl}) are
smooth functions of local coordinates,
$K$ is an arbitrary constant, were introduced and studied.
These nonlocal Poisson brackets play an important
role in the theory of systems of hydrodynamic type.
The form of bracket (\ref{nonl}) is invariant
with respect to local changes of coordinates.
A bracket of form (\ref{nonl}) is said to be
{\it nondegenerate} if
$\det (g^{ij} (u)) \not\equiv 0$.
If $\det (g^{ij} (u)) \not\equiv 0$,
then bracket (\ref{nonl}) is a Poisson bracket if and only if
$g^{ij} (u)$ is an arbitrary pseudo-Riemannian
contravariant metric of constant Riemannian curvature
$K$,
$b^{ij}_k (u) = - g^{is} (u) \Gamma ^j_{sk} (u),$ where
$\Gamma^j_{sk} (u)$ is the Riemannian connection
generated by the metric $g^{ij} (u)$
(the Levi--Civita connection) \cite{26x}
(note that the coefficients
$g^{ij} (u)$ and $b^{ij}_k (u)$ of bracket
(\ref{nonl}) are transformed as corresponding
differential-geometric objects under local changes of
coordinates: a contravariant metric
$g^{ij} (u)$ and a connection $b^{ij}_k (u)
= - g^{is} (u) \Gamma^j_{sk} (u)$, respectively, $K$
is an invariant).
For $K= 0$ we have the local Poisson brackets
of hydrodynamic type
(the Dubrovin--Novikov brackets \cite{25x}).

The description problem for compatible metrics of
constant Riemannian curvature is equivalent to
that for compatible nonlocal Poisson brackets
of hydrodynamic type generated by metrics of
constant Riemannian curvature (compatible Mokhov--Ferapontov
brackets).
Recall that Poisson brackets are said to be
{\it compatible} if their linear combination is also
a Poisson bracket
(Magri, \cite{27x}).
As was shown in \cite{28x}, \cite{29x}
(see also \cite{16x}, \cite{30x}, \cite{31x}),
compatible Mokhov--Ferapontov brackets are described by
a consistent nonlinear system integrable by
the method of the inverse scattering problem
(the case of compatible Dubrovin--Novikov
brackets, i.e., compatible local Poisson brackets of
hydrodynamic type, was integrated by the author earlier in
\cite{13x}, \cite{14x}, see also \cite{16x},
\cite{21x}).

\section{Compatible nonlocal Poisson brackets
 \\ of hydrodynamic type}

\bl \label{locnon}
In the classification problem for an arbitrary
pair of compatible nonlocal Poisson brackets of form
(\ref{nonl}), one can always consider one of these
two Poisson brackets as local without loss of generality.
\el

Actually, if two compatible nonlocal Poisson brackets
of form (\ref{nonl})
$\{ I, J \}_0$ (with a corresponding constant $K_0$
in the nonlocal term) and $\{ I, J \}_1$ (with a constant $K_1$)
are linear independent, then in the pencil of these Poisson brackets,
i.e., among the Poisson brackets $\lambda_0 \{ I, J \}_0 +
\lambda_1 \{ I, J \}_1$, where $\lambda_0$ and $\lambda_1$
are arbitrary constants, there is necessarily a nonzero
local Poisson bracket
$ \{ I, J \} = \lambda'_0 \{ I, J \}_0 + \lambda'_1 \{ I, J \}_1,$
where $ \lambda'_0$ and $\lambda'_1$ are arbitrary constants
satisfying the relation $\lambda_0 K_0 + \lambda_1 K_1 = 0$,
which can be taken as one of the generators
for all the considered pencil of compatible
Poisson brackets (it is obvious that if
this local Poisson bracket $\{ I, J \}$ is identically zero,
then the Poisson brackets $ \{ I, J \}_0$ and $ \{ I, J \}_1$
are linearly dependent: $ \lambda'_0 \{ I, J \}_0 +
\lambda'_1 \{ I, J \}_1 \equiv 0$).

Consider the problem of compatibility
for a pair of nonlocal and local Poisson
brackets of hydrodynamiic type
\be
\{ I,J \}_1 = \int {\delta I \over \delta u^i(x)}
\left ( g^{ij}_1 (u(x)) {d \over dx} + b^{ij}_{1, k} (u(x))\, u^k_x
+ K_1 u^i_x \left ( {d \over dx} \right )^{-1} u^j_x \right )
{\delta J \over \delta u^j(x)} dx
\label{nonloc}
\ee
and
\be
\{ I,J \}_2 = \int {\delta I \over \delta u^i(x)}
\left ( g^{ij}_2 (u(x)) {d \over dx} + b^{ij}_{2, k} (u(x))\, u^k_x
 \right )
{\delta J \over \delta u^j(x)} dx,
\label{lok}
\ee
i.e., the condition that for any constant $\lambda$
the bracket
\be
\{ I,J \} = \{ I, J \}_1 + \lambda \{ I, J \}_2
\label{comp}
\ee
is a Poisson bracket (thus, formula (\ref{comp})
defines {\it a pencil of compatible Poisson brackets}).

We assume further that the local
bracket $\{ I, J \}_2$ is nondegenerate, i.e.,
$\det (g^{ij}_2 (u)) \not\equiv 0$,
but we do not impose any additional conditions
on the bracket $\{ I, J \}_1$, i.e.,
generally speaking, this bracket may be degenerate.
Bracket (\ref{comp})
can be degenerate, therefore
here we indicate the general relations
for the coefficients of a bracket of form (\ref{nonl})
which are equivalent to the condition
that bracket (\ref{nonl}) is a Poisson bracket.
These general relations (without
the assumption of nondegeneracy)
were obtained in the present author's work
\cite{32x}
(see also \cite{33x}, \cite{34x}):
\be
g^{ij} (u) = g^{ji} (u), \label{s1}
\ee
\be
{\pa g^{ij} \over \pa u^k} = b^{ij}_k (u) + b^{ji}_k (u),\label{s2}
\ee
\be
g^{is} (u) b^{jr}_s (u) = g^{js} (u) b^{ir}_s (u), \label{s3}
\ee
\be
g^{is} (u) \left ( {\pa b^{jr}_s \over \pa u^k}  -
{\pa b^{jr}_k \over \pa u^s} \right )
+ b^{ij}_s (u) b^{sr}_k (u) - b^{ir}_s (u) b^{sj}_k (u) =
K ( g^{ir} (u) \delta^j_k -
g^{ij} (u) \delta^r_k ),\label{s4}
\ee
\bea
&&
\sum_{(i, j, r)} \left [
b^{si}_p (u) \left (
{\pa b^{jr}_k \over \pa u^s} -
{\pa b^{jr}_s \over \pa u^k} \right ) + K (b^{ij}_k (u) - b^{ji}_k (u))
\delta^r_p +\right. \nn\\
&&
\left. b^{si}_k (u) \left (
{\pa b^{jr}_p \over \pa u^s} -
{\pa b^{jr}_s \over \pa u^p} \right ) + K (b^{ij}_p (u) - b^{ji}_p (u))
\delta^r_k \right ]  = 0, \label{s5}
\eea
where $\sum_{(i, j, r)}$ means summation over all
cyclic permutations of the indices $i, j, r$.

\section{Canonical form for compatible pairs of
brackets}

According to the Dubrovin--Novikov theorem \cite{25x},
for any nondegenerate local Poisson bracket
of hydrodynamic type $\{ I, J \}_2$, there always
exit local coordinates $u^1,...,u^N$
(flat coordinates of the metric $g^{ij}_2 (u)$)
in which this bracket is constant, i.e.,
$g^{ij}_2 (u) = \eta^{ij} = {\rm\ const}$, $b^{ij}_{2, k} (u) =
\Gamma^i_{2, jk} (u) = 0$.
Thus we can choose flat coordinates of the metric
$g^{ij}_2 (u)$ and further we consider that
the Poisson bracket
$\{ I, J \}_2$ is constant and has the form
\be
\{ I,J \}_2 = \int {\delta I \over \delta u^i(x)}
 \eta^{ij} {d \over dx}
{\delta J \over \delta u^j(x)} dx,
\label{lokconst}
\ee
where $\eta^{ij} = \eta^{ji}$, $\eta^{ij} = {\rm\ const},$
$\det (\eta^{ij}) \neq 0$.
In the sequel, in the considered flat coordinates,
we use also the covariant metric $\eta_{ij}$,
which is inverse to the contravariant
metric $\eta^{ij}$: $\eta^{is} \eta_{sj} = \delta^i_j.$

\bt [\cite{30x}] \label{t1}
An arbitrary nonlocal Poisson bracket $\{ I, J \}_1$
of form (\ref{nonloc})
(may be degenerate) is compatible with the constant
Poisson bracket
(\ref{lokconst}) if and only if it has
the form
\bea
&&
\ \{ I,J \}_1 = \int {\delta I \over \delta u^i(x)}
\left ( \left [ \eta^{is} {\pa H^j \over
\pa u^s} + \eta^{js} {\pa H^i \over \pa u^s}
- K_1 u^i u^j \right ] {d \over dx} + \right.\nn\\
&&
\left. \left [
\eta^{is} {\pa^2 H^j \over \pa u^s \pa u^k} - K_1 \delta^i_k u^j
\right ]\, u^k_x
+ K_1 u^i_x \left ( {d \over dx} \right )^{-1} u^j_x \right )
{\delta J \over \delta u^j(x)} dx,
\label{nonloc3}
\eea
where $H^i (u),$ $1 \leq i \leq N,$ are smooth functions
defined in a certain domain of local coordinates.
\et

In the flat case of compatible Dubrovin--Novikov
brackets ($K_1 = 0$), the corresponding statement
was formulated and proved by the present author
in \cite{35x}--\cite{37x}
(see also the conditions on flat pencils of metrics in \cite{9x}).

\section{Integrable equations for canonical
\\ compatible pairs of brackets}

\bt [\cite{30x}]  \label{liu}
An arbitrary nonlocal bracket of form (\ref{nonloc3})
(may be degenerate) is a Poisson bracket if and only if
the following equations are satisfied:
\be
{\pa^2 H^i \over \pa u^k \pa u^s} \eta^{sp}
{\pa^2 H^j \over \pa u^p \pa u^l} =
{\pa^2 H^j \over \pa u^k \pa u^s} \eta^{sp}
{\pa^2 H^i \over \pa u^p \pa u^l}, \label{ass1}
\ee
\bea
&&
\left ( \eta^{ir} {\pa H^s \over \pa u^r} +
\eta^{sr} {\pa H^i \over \pa u^r} - K_1 u^i u^s \right ) \eta^{jp}
{\pa^2 H^k \over \pa u^p \pa u^s} =\nn\\
&&
\left ( \eta^{jr} {\pa H^s \over \pa u^r} +
\eta^{sr} {\pa H^j \over \pa u^r} - K_1 u^j u^s
 \right ) \eta^{ip}
{\pa^2 H^k \over \pa u^p \pa u^s}. \label{ass2}
\eea
\et

In the flat case ($K_1 = 0$), the corresponding theorem
was obtained by the present author in \cite{36x},
where was also stated the conjecture on the integrability
of system (\ref{ass1}),
(\ref{ass2}) for $K_1 = 0$ by the method of the inverse scattering
problem. This conjecture was proved by the author in the works
\cite{13x}, \cite{14x},
\cite{16x} (see also
\cite{21x}, where a Lax pair for system (\ref{ass1}),
(\ref{ass2}) was indicated). The corresponding
general conditions on flat pencils of metrics
were indicated in \cite{9x}.

It is obvious that any set of $N$ linear functions
$H^i (u) = a^i_k u^k + a^i,$ where $a^i_k$ is an
arbitrary constant matrix,
$a^i_k = const,$ $a^i = const,$ is always a solution
of the nonlinear system (\ref{ass1}),
(\ref{ass2}) and, consequently, always generates
a corresponding canonical pair of compatible
Poisson brackets
(\ref{lokconst}), (\ref{nonloc3}) (see \cite{30x}).

In the flat case (for $K_1 = 0$), a set of $N$
quadratic functions $H^i (u) = c^i_{jk} u^j u^k,$
where $c^i_{jk} = c^i_{kj},$ $c^i_{jk} = {\rm const},$
is a solution of nonlinear system (\ref{ass1}),
(\ref{ass2}) if and only if the structural constants
$a^{ij}_k = \eta^{is} c^j_{sk}$ satisfy the relations (see \cite{10x})
\be
a^{ij}_s a^{sk}_l = a^{ik}_s a^{sj}_l,
\ee
\be
(a^{is}_l + a^{si}_l) a^{jk}_s =
(a^{js}_l + a^{sj}_l) a^{ik}_s,
\ee
which are equivalent to the condition that
$N$-dimensional algebra with basis $e^1,...,e^N$
and multiplication
$$e^i \cdot e^j = a^{ij}_k e^k$$ is a Novikov algebra, i.e.,
the identity
\be
(e^i \cdot e^j) \cdot e^k =
(e^i \cdot e^k) \cdot e^j,
\ee
\be
e^i \cdot (e^j \cdot e^k) - (e^i \cdot e^j) \cdot e^k =
e^j \cdot (e^i \cdot e^k) - (e^j \cdot e^i) \cdot e^k  \label{levo}
\ee
are fulfilled (see \cite{10x}).
Moreover,
the invariant nondegenerate symmetric bilinear form
$$( e^i, e^j ) = \eta^{ij}$$
is defined on this Novikov algebra, i.e.,
$$( e^i \cdot e^j, e^k ) =
( e^i, e^k \cdot e^j ), \ \ \ ( e^i, e^j ) = ( e^j, e^i ).$$
Thus, in this case, we get a special class
of quasi-Frobenius algebras, namely, exactly
the class of left-symmetric quasi-Frobenius algebras
describing linear one-dimensional local Poisson brackets of
hydrodynamic type (see \cite{4x}).
\bl
The condition of left-symmetry (\ref{levo})
in an arbitrary algebra
${\cal {F}}$ with multiplication
$e^i \cdot e^j = f^{ij}_k e^k$ and identity (\ref{q1})
(and consequently in any quasi-Frobenius
algebra) is equivalent to the condition that
the symmetric bilinear form
\be
< e^i, e^j > = (f^{ij}_k + f^{ji}_k) u^k
\ee
is invariant for any $u = (u^1,..., u^N)$, i.e.,
\be
< e^i \cdot e^j, e^k > =
< e^i, e^k \cdot e^j >.
\ee
The parameters $u^1,..., u^N$ realize a deformation of
the invariant form.
\el

In the general case of compatible metrics
of constant Riemannian curvature,
it is interesting to study algebraic structures
related to sets of $N$ cubic functions $H^i (u) =
c^i_{jkl} u^j u^k u^l,$ where $c^i_{jkl} = c^i_{kjl} =
c^i_{jlk},$ $c^i_{jkl} = {\rm const}.$
Such a set of cubic functions is a solution of
nonlinear system (\ref{ass1}),
(\ref{ass2}) if and only if for the structural
constants $a^{ij}_{kl} = \eta^{is} c^j_{skl},$
$a^{ij}_{kl} = a^{ij}_{lk},$ the following relations
are satisfied:
\be
a^{ki}_{ms} a^{sj}_{nl} - a^{kj}_{ns} a^{si}_{ml}
+ a^{ki}_{ns} a^{sj}_{ml} - a^{kj}_{ms} a^{si}_{nl} = 0,
\ee
\be
\sum_{[l, m, n]} [ (3 a^{is}_{mn} + 3 a^{si}_{mn} -
K_1 \delta^i_m \delta^s_n) a^{jk}_{ls} -
(3 a^{js}_{mn} + 3 a^{sj}_{mn} -
K_1 \delta^j_m \delta^s_n) a^{ik}_{ls} ] = 0,
\ee
where $\sum_{[l, m, n]}$ means summation over all
permutations of the indices $l, m, n.$
In this case, the identity
\be
(a \ast b) \ast c =
(a \ast c) \ast b
\ee
is fulfilled in $N$-dimensional algebra ${\cal {C}} (u)$
with basis $e^1,...,e^N$ and multiplication
$e^i \ast e^j = a^{ij}_{lk} u^l e^k$ for any
$u = (u^1,..., u^N)$, and, moreover,
this algebra ${\cal {C}} (u)$ possesses two
invariant symmetric bilinear forms
$( e^i, e^j ) = \eta^{ij}$ and
\be
< e^i, e^j > = 3 (a^{ij}_{kl} + a^{ji}_{kl}) u^k u^l -
 K_1 u^i u^j,
\ee
i.e.,
$$( e^i \ast e^j, e^k ) =
( e^i, e^k \ast e^j ), \ \ \
< e^i \ast e^j, e^k > =
< e^i, e^k \ast e^j >.$$
Thus, in this case, we get
$N$-parametric deformations of quasi-Frobenius
algebras
$({\cal {C}} (u), ( \cdot , \cdot ))$ and
$({\cal {C}} (u), < \cdot , \cdot >).$

\bt [\cite{28x}, \cite{29x}] \label{int}
The system of nonlinear equations (\ref{ass1}), (\ref{ass2})
is consistent and integrable by the method of
the inverse scattering problem.
\et

Note that the system
of nonlinear equations that was found and integrated
in \cite{28x}, \cite{29x} is equivalent to system
(\ref{ass1}), (\ref{ass2}) and, consequently,
also describes compatible nonlocal Poisson brackets of
form (\ref{nonl}), but in different special local
coordinates, which are much more convenient
for the integration (the metrics of both compatible
brackets are diagonal in these coordinates).
In the flat case ($K_1 = 0$), the corresponding
system in special ``diagonal'' local coordinates was
integrated by the author in \cite{13x}, \cite{14x}.

\section{Quasi-Frobenius algebras and \\ their
integrable deformations}

Consider the Poisson bracket (\ref{nonloc3})
defining the canonical pair of compatible brackets and
for any $u = (u^1,..., u^N)$  define an algebra
${\cal {A}} (u)$ in $N$-dimensional vector space
with basis
$e^1,..., e^N$ and multiplication
\be
e^i \circ e^j = \eta^{is} {\pa^2 H^j \over \pa u^s \pa u^k } \, e^k.
\ee
Moreover, define a symmetric bilinear form
(``a metric of constant Riemannian curvature
$K_1$'')
\be
< e^i, e^j > = g^{ij}_1 (u)
\ee
on the algebra ${\cal {A}} (u)$, i.e.,
\be
< e^i, e^j > = \eta^{is} {\pa H^j \over \pa u^s}
+ \eta^{js} {\pa H^i \over \pa u^s} - K_1 u^i u^j.
\ee
Then equations (\ref{ass1}) mean that the identity
\be
(e^i \circ e^j) \circ e^k =
(e^i \circ e^k) \circ e^j  \label{soot1}
\ee
is fulfilled in the algebra ${\cal {A}} (u)$,
and equations (\ref{ass2}) are equivalent to the identity
\be
< e^i \circ e^j, e^k > = < e^i, e^k \circ e^j >. \label{soot2}
\ee
Thus, any solution of the integrable nonlinear system
(\ref{ass1}), (\ref{ass2}) defines a quasi-Frobenius algebra
${\cal {A}} (u)$ with the multiplication $a \circ b$
and the invariant symmetric bilinear form
$< a, b >$ for any fixed $u = (u^1,...,u^N)$,
and system (\ref{ass1}), (\ref{ass2}) defines
an $N$-parametric deformation of this quasi-Frobenius algebra.
Note that also there is always a nondeformed
invariant nondegenerate symmetric bilinear form
(``a flat metric'')
\be
( e^i, e^j ) = \eta^{ij}, \ \ \
(e^i \circ e^j, e^k) = (e^i, e^k \circ e^j)
\ee
on the constructed quasi-Frobenius algebra ${\cal {A}} (u)$.
Thus, the compatible metrics of constant Riemannian
curvature $\eta^{ij}$ and
$\eta^{is} {\pa H^j \over \pa u^s}
+ \eta^{js} {\pa H^i \over \pa u^s} - K_1 u^i u^j$
define two ``compatible'' quasi-Frobenius algebras
$({\cal {A}} (u), ( \cdot , \cdot ))$ and
$({\cal {A}} (u), < \cdot , \cdot >)$
for any $u = (u^1,..., u^N)$.
Quasi-Frobenius structures and quasi-Frobenius
manifolds arising in the flat case for
$K_1 = 0$ were considered in \cite{10x}.

{\bf Conjecture}. The deformations of noncommutative
associative algebras, recently found by S.M.Natanzon in
 \cite{22x} (see also \cite{23x})
and corresponding to open-closed two-dimensional
topological field theories, belong to an integrable
class of deformations of noncommutative
quasi-Frobenius algebras and are integrable by
the method of the inverse scattering problem.

The construction of the corresponding integrable ``associative''
reductions is a very interesting and important problem.

\section{Deformations of Frobenius algebras}

Assume that multiplication in the algebra
${\cal {A}} (u)$ is commutative:
$e^i \circ e^j = e^j \circ e^i,$ i.e.,
\be
\eta^{is} {\pa^2 H^j \over \pa u^s \pa u^k} =
\eta^{js} {\pa^2 H^i \over \pa u^s \pa u^k}.
\ee

Then there exist $c_{pl} = {\rm const}$ such that
\be
\eta_{ps} {\pa H^s \over \pa u^l} =
\eta_{ls} {\pa H^s \over \pa u^p} + c_{lp} - c_{pl},
\ee
and, consequently,
\be
{\pa (\eta_{ps} H^s + c_{ps} u^s) \over \pa u^l} =
{\pa (\eta_{ls} H^s + c_{ls} u^s) \over \pa u^p}.
\ee
Thus, there exists a function $\Phi (u)$
(``a potential'') such that
\be
H^i (u) = \eta^{is} \left ( {\pa \Phi \over \pa u^s} -
{1 \over 2} (c_{sk} - c_{ks}) u^k \right ), \label{re1}
\ee
where we use that the symmetric part of the matrix
$(c_{ks})$ can easily be included in
``the potential'' $\Phi (u)$:
$$ {1 \over 2} (c_{sk} + c_{ks}) u^k = {1 \over 2}
{\partial (c_{pr} u^p u^r) \over \partial u^s}.$$

Thus, for the commutative algebra ${\cal {A}} (u)$,
we have:
\be
e^i \circ e^j = \eta^{is} \eta^{jp} {\pa^3 \Phi \over
\pa u^s \pa u^p \pa u^k} \, e^k,
\ee
\be
< e^i, e^j > = 2 \eta^{is} \eta^{jp} {\pa^2 \Phi \over
\pa u^s \pa u^p} - K_1 u^i u^j.
\ee

It follows from identities
(\ref{soot1}) and (\ref{soot2}) that
the commutative algebra ${\cal {A}} (u)$ is associative:
\be
(e^i \circ e^j) \circ e^k = e^i \circ ( e^j \circ e^k),
\ee
and is equipped with an invariant symmetric
bilinear form:
\be
< e^i \circ e^j, e^k > = < e^i, e^j \circ e^k >,
\ee
i.e., for any $u = (u^1,...,u^N)$
the algebra ${\cal {A}} (u)$ is Frobenius.

As a result of the reduction (\ref{re1}),
equations (\ref{ass1}) assume the form of
the Witten--Dijkgraaf--Verlinde--Verlinde--Dubrovin
equations of associativity
(see \cite{9x}, \cite{38x}--\cite{41x}):
\be
{\pa^3 \Phi \over \pa u^k \pa u^i \pa u^s} \eta^{sp}
{\pa^3 \Phi \over \pa u^p \pa u^j \pa u^l} =
{\pa^3 \Phi \over \pa u^k \pa u^j \pa u^s} \eta^{sp}
{\pa^3 \Phi \over \pa u^p \pa u^i \pa u^l}, \label{as1}
\ee
and equations (\ref{ass2}) have the following form in this case:
\bea
&&
\left ( {\pa^2 \Phi \over \pa u^i \pa u^s} -
 {K_1 \over 2} \eta_{ir} \eta_{sl} u^r u^l \right ) \, \eta^{sp} \,
{\pa^3 \Phi \over \pa u^p \pa u^j \pa u^k} =\nn\\
&&
\left ( {\pa^2 \Phi \over \pa u^j \pa u^s} -
 {K_1 \over 2} \eta_{jr} \eta_{sl} u^r u^l \right ) \, \eta^{sp} \,
{\pa^3 \Phi \over \pa u^p \pa u^i \pa u^k}.
\label{as2}
\eea

Note that in the flat case (for $K_1 = 0$)
all Dubrovin's potentials
$\Phi (u)$ corresponding to two-dimensional topological
field theories are always solutions for both
corresponding equations (\ref{as1}) and (\ref{as2}) with $K_1 =0$
(see \cite{9x}).

\begin{flushleft}
Centre for Nonlinear Studies,\\
L.D.Landau Institute for Theoretical Physics,\\
Russian Academy of Sciences\\
e-mail: mokhov@mi.ras.ru, mokhov@landau.ac.ru\\
\end{flushleft}


\begin{thebibliography}{99}



\bibitem{1x} Gelfand I.M., Dorfman I.Ya.
Hamiltonian operators and algebraic structures associated
with them. Funkts. Anal. Pril. 1979. Vol. 13, No. 4.
P. 13--30; English transl. in Funct. Anal. Appl.
1979. Vol. 13, No. 4. P. 246--262.


\bibitem{2x} Gelfand I.M., Dorfman I.Ya.
Hamiltonian operators and infinite-dimensional Lie algebras.
Funkts. Anal. Pril. 1981. Vol. 15, No. 3.
P. 23--40; English transl. in Funct. Anal. Appl.
1981. Vol. 15, No. 3. P. 173--187.


\bibitem{3x} Dorfman I. Dirac structures and integrability
of nonlinear evolution equations. Chichester, England. Wiley (1993).


\bibitem{4x} Balinskii A.A., Novikov S.P.
Poisson brackets of hydrodynamic type,
Frobenius algebras and Lie algebras.
Dokl. Akad. Nauk SSSR. 1985. Vol. 283, No. 5.
P. 1036--1039; English transl. in Soviet Math. Dokl.
1985. Vol. 32. P. 228--231.



\bibitem{5x} Novikov S.P. The geometry of conservative
systems of hydrodynamic type. The method of averaging
for field-theoretical systems. Uspekhi Mat. Nauk. 1985.
Vol. 40, No. 4. P. 79--89; English transl. in
Russian Math. Surveys. 1985. Vol. 40, No. 4. P. 85--98.


\bibitem{6x} Zelmanov E.I. On a class of local
translation-invariant Lie algebras.
Dokl. Akad. Nauk SSSR. 1987. Vol. 292, No. 6.
P. 1294--1297; English transl. in Soviet Math. Dokl.
1987. Vol. 35. P. 216--218.


\bibitem{7x} Vinberg E.B. Convex homogeneous domains.
Dokl. Akad. Nauk SSSR. 1961. Vol. 141, No. 3.
P. 521--524; English transl. in Soviet Math. Dokl.
1961.


\bibitem{8x} Mokhov O.I. On Poisson brackets of
Dubrovin--Novikov type (DN-brackets).
Funkts. Anal. Pril. 1988. Vol. 22, No. 4.
P. 92--93; English transl. in Funct. Anal. Appl.
1988. Vol. 22, No. 4. P. 336--338.




\bibitem{9x} Dubrovin B. Geometry of 2D topological field theories.
Lecture Notes in Math. 1996. Vol. 1620. P. 120--348;
E-print. arXiv: hep-th/9407018 (1994).



\bibitem{10x} Dubrovin B., Zhang Y. Normal forms
of hierarchies of integrable PDEs,
Frobenius manifolds and Gromov--Witten invariants.
E-print. arXiv: math.DG/0108160 (2001).



\bibitem{11x} Mokhov O.I.
Compatible Dubrovin--Novikov Hamiltonian operators,
Lie derivative and integrable systems of hydrodynamic type.
Proceedings of the International Conference ``Nonlinear
Evolution Equations and Dynamical Systems.'' Cambridge (England), July
24--30, 2001 (will be published in Theoretical and Mathematical
Physics, 2002);
E-print. arXiv: math.DG/0201281 (2002).



\bibitem{12x} Mokhov O.I. Compatible Dubrovin--Novikov
Hamiltonian operators and the Lie derivative.
Uspekhi Mat. Nauk. 2001.
Vol. 56, No. 6. P. 161--162; English transl. in
Russian Math. Surveys. 2001. Vol. 56, No. 6. P. 1175--1176.



\bibitem{12xx} Dubrovin B. Flat pencils of metrics and
Frobenius manifolds. Preprint SISSA 25/98/FM;
E-print. arXiv: math.DG/9803106 (1998).





\bibitem{13x} Mokhov O.I. Integrability of the equations
for nonsingular pairs of compatible flat metrics.
Teoretich. i Matematich. Fizika. 2002. Vol. 130, No. 2.
P. 233--250; English transl. in Theoretical and Mathematical
Physics. 2002. Vol. 130, No. 2. P. 198--212; E-print.
arXiv: math.DG/0005081 (2000).


\bibitem{14x} Mokhov O.I. Flat pencils of metrics
and integrable reductions of the Lam\'e
equations. Uspekhi Mat. Nauk. 2001.
Vol. 56, No. 2. P. 221--222; English transl. in
Russian Math. Surveys. 2001. Vol. 56, No. 2.



\bibitem{15x} Mokhov O.I.
Integrable bi-Hamiltonian  systems of hydrodynamic type.
Uspekhi Mat. Nauk. 2002.
Vol. 57, No. 1. P. 157--158; English transl. in
Russian Math. Surveys. 2002. Vol. 57, No. 1.


\bibitem{16x} Mokhov O.I.
Compatible and almost compatible pseudo-Riemannian
metrics.
Funkts. Anal. Pril. 2001. Vol. 35, No. 2.
P. 24--36; English transl. in Funct. Anal. Appl.
2001. Vol. 35, No. 2. P. 100--110;
E-print. arXiv: math.DG/0005051 (2000).



\bibitem{21x} Ferapontov E.V. Compatible Poisson brackets
of hydrodynamic type. E-print. arXiv: math.DG/0005221 (2000).


\bibitem{22x} Natanzon S.M. Extension
cohomological fields theory and
noncommutative Frobenius manifolds.
E-print. arXiv: math-ph/0206033 (2002).


\bibitem{23x} Alexeevski A., Natanzon S.
Non-commutative extensions of
two-dimensional topological field theories
and Hurwitz numbers for real algebraic curves.
E-print. arXiv: math.GT/0202164 (2002).




\bibitem{24x} Mokhov O.I.
Compatible and almost compatible metrics.
Uspekhi Mat. Nauk. 2000.
Vol. 55, No. 4. P. 217--218; English transl. in
Russian Math. Surveys. 2000. Vol. 55, No. 4. P. 819--821.


\bibitem{25x} Dubrovin B.A., Novikov S.P.
The Hamiltonian formalism of one-dimensional
systems of hydrodynamic type and the
Bogolyubov--Whitham averaging method.
Dokl. Akad. Nauk SSSR. 1983. Vol. 270, No. 4.
P. 781--785; English transl. in Soviet Math. Dokl.
1983. Vol. 27. P. 665--669.



\bibitem{26x} Mokhov O.I., Ferapontov E.V.
Non-local Hamiltonian operators of hydrodynamic type
related to metrics of constant curvature.
Uspekhi Mat. Nauk. 1990.
Vol. 45, No. 3. P. 191--192; English transl. in
Russian Math. Surveys. 1990. Vol. 45, No. 3. P. 218--219.




\bibitem{27x}  Magri F. A simple model of the integrable
Hamiltonian equation. J. Math. Phys. 1978.
Vol. 19, No. 5. P. 1156--1162.



\bibitem{28x} Mokhov O.I.
Compatible metrics of constant Riemannian curvature:
local geometry, nonlinear equations, and integrability.
Funkts. Anal. Pril. 2002. Vol. 36, No. 3.
P. 36--47; English transl. in Funct. Anal. Appl.
2002. Vol. 36, No. 3;
E-print. arXiv: math.DG/0201280 (2002).


\bibitem{29x} Mokhov O.I.
Lax pair for nonsingular pencils of metrics of
constant Riemannian curvature.
Uspekhi Mat. Nauk. 2002.
Vol. 57, No. 3. P. 155--156; English transl. in
Russian Math. Surveys. 2002. Vol. 57, No. 3.


\bibitem{30x} Mokhov O.I.
Liouville canonical form for compatible
nonlocal Poisson brackets of hydrodynamic type,
and integrable hierarchies.
E-print. arXiv: math.DG/0201223 (2002).


\bibitem{31x} Mokhov O.I.
Integrable bi-Hamiltonian hierarchies
generated by compatible metrics of constant
Riemannian curvature.
Uspekhi Mat. Nauk. 2002.
Vol. 57, No. 5; English transl. in
Russian Math. Surveys. 2002. Vol. 57, No. 5.











\bibitem{32x}  Mokhov O.I. Hamiltonian systems of
hydrodynamic type and constant curvature metrics.
Phys. Letters A. 1992. Vol. 166, Nos. 3-4. P. 215--216.



\bibitem{33x} Mokhov O.I.
Symplectic and Poisson structures on loop
spaces of smooth manifolds, and integrable systems.
Uspekhi Mat. Nauk. 1998.
Vol. 53, No. 3. P. 85--192; English transl. in
Russian Math. Surveys. 1998. Vol. 53, No. 3. P. 515--622.



\bibitem{34x} Mokhov O.I. Symplectic and Poisson Geometry
on Loop Spaces of Smooth Manifolds and Integrable Equations.
Reviews in Mathematics and Mathematical Physics.
Eds. S.P.Novikov and I.M.Krichever. 2001. Vol. 11, part 2.
Harwood Academic Publishers. P. 1--204.


\bibitem{35x} Mokhov O.I.
On compatible Poisson structures of
hydrodynamic type.
Uspekhi Mat. Nauk. 1997.
Vol. 52, No. 6. P. 171--172; English transl. in
Russian Math. Surveys. 1997. Vol. 52, No. 6. P. 1310--1311.



\bibitem{36x} Mokhov O.I.
Compatible Poisson structures
of hydrodynamic type and associativity equations.
Trudy Matem. Inst. Akad. Nauk. 1999. Vol. 225.
Moscow. Nauka. P. 284--300;
English transl. in
Proceedings of the Steklov Institute of
Mathematics (Moscow). 1999. Vol. 225. P. 269--284.




\bibitem{37x}  Mokhov O.I. Compatible Poisson structures
of hydrodynamic type and the equations of
associativity in two-dimensional topological field theory.
Reports on Mathematical Physics. 1999. Vol. 43, No. 1/2.
P. 247--256.



\bibitem{38x} Witten E. On the structure of topological phase
of two-dimensional gravity. Nucl. Phys. B. 1990. Vol. 340.
P. 281--332.


\bibitem{39x} Witten E. Two-dimensional gravity and
intersection theory on moduli space. Surveys in Diff. Geometry.
1991. Vol. 1. P. 243--310.


\bibitem{40x} Dijkgraaf R., Verlinde H., Verlinde E.
Topological strings in $d < 1$. Nucl. Phys. B. 1991.
Vol. 352. P. 59--86.


\bibitem{41x} Dubrovin B. Integrable systems in topological
field theory. Nucl. Phys. B. 1992. Vol. 379. P. 627--689.


\end{thebibliography}
\end{document}